# OPTIMAL CHANGE-POINT ESTIMATION FROM INDIRECT OBSERVATIONS

BY A. GOLDENSHLUGER,[1] A. TSYBAKOV AND A. ZEEVI[3]

*University of Haifa, Université Paris VI and Columbia University*

*Dedicated to Boris Polyak on the occasion of his 70th birthday*

We study nonparametric change-point estimation from indirect noisy observations. Focusing on the white noise convolution model, we consider two classes of functions that are smooth apart from the change-point. We establish lower bounds on the minimax risk in estimating the change-point and develop rate optimal estimation procedures. The results demonstrate that the best achievable rates of convergence are determined both by smoothness of the function away from the change-point and by the degree of ill-posedness of the convolution operator. Optimality is obtained by introducing a new technique that involves, as a key element, detection of zero crossings of an estimate of the properly smoothed second derivative of the underlying function.

**1. Introduction.** In this paper we study the problem of change-point estimation from indirect and noisy observations. Let $f \in \mathbb{L}_2(\mathbb{R})$ denote the unknown function. Consider the white noise model

$$(1) \qquad dY(x) = (\mathsf{K}f)(x)\,dx + \varepsilon\,dW(x), \qquad x \in \mathbb{R},$$

where $W(\cdot)$ is the standard two-sided Wiener process on $\mathbb{R}$, $0 < \varepsilon < 1$, and $\mathsf{K}$ is the convolution operator with kernel $K \in \mathbb{L}_1(\mathbb{R})$ whose action on a function $f \in \mathbb{L}_2(\mathbb{R})$ is defined by

$$(2) \qquad (\mathsf{K}f)(x) = \int_{-\infty}^{\infty} K(x-y)f(y)\,dy.$$

Received July 2004; revised February 2005.
[1]Supported by Israel Science Foundation Grant 300/04.
[3]Supported in part by NSF Grant 04-47652.

We assume that $f$ is smooth apart from a jump discontinuity of the first kind at a point $\theta$ and, without loss of generality, we suppose that $\theta \in [0, 1]$. The problem is to estimate the change-point $\theta$ based on the observation of a trajectory of the process $Y(\cdot)$ that satisfies (1).

We study this problem in a minimax framework. Let $\hat{\theta}$ be an estimator of $\theta$ based on observation of $Y(\cdot)$ that satisfies (1). We measure the accuracy of $\hat{\theta}$ by the maximal risk

$$R_\varepsilon[\hat{\theta}; \mathcal{G}] = \sup_{f \in \mathcal{G}} \{\mathbb{E}_f |\hat{\theta} - \theta|^2\}^{1/2}$$

over a class of functions $\mathcal{G}$ that have a single change-point $\theta \in [0, 1]$. Here $\mathbb{E}_f$ denotes the expectation with respect to the probability distribution $\mathbb{P}_f$ generated by the model (1) with given $f$. The minimax risk is defined by

$$R_\varepsilon^*[\mathcal{G}] = \inf_{\hat{\theta}} R_\varepsilon[\hat{\theta}; \mathcal{G}],$$

where the infimum is taken over all possible estimators of $\theta$. An estimator $\tilde{\theta}$ is called *rate optimal* on the class $\mathcal{G}$ if it satisfies

$$R_\varepsilon[\tilde{\theta}; \mathcal{G}] \asymp R_\varepsilon^*[\mathcal{G}] \qquad \text{as } \varepsilon \to 0.$$

Our aim is to find rate optimal change-point estimators and to establish asymptotics of minimax risks for some natural classes of functions $\mathcal{G}$ and operators $\mathsf{K}$.

Change-points and singularities are intrinsic features of signals that appear in a wide variety of applied contexts in economics, medicine and physical science. For many types of signals, change-points convey important information about underlying phenomena. For instance, in images, discontinuities of the intensity function correspond to the location of the contour of an object that may be particularly important for recognition purposes. We refer to the volume by Carlstein, Müller and Siegmund [4] for a comprehensive survey of the area and references. The problem of nonparametric change-point estimation has been extensively studied in the case where the observations of $f$ are *direct*, that is, when $\mathsf{K}$ is the identity operator. For such a model, Korostelev [15] constructed a rate optimal estimator of $\theta$ and derived the optimal rates of convergence. He showed that the minimax risk over the class of functions that have a single change-point and satisfy the Lipschitz condition away from the change-point, converges to zero at the rate $\varepsilon^2$, which is faster than the usual parametric rate. (Here and in what follows we have in mind a standard correspondence between the Gaussian white noise model and discrete sample models (cf. [3]), given by the calibration $\varepsilon = n^{-1/2}$, where $n$ is the sample size. With this calibration, the term parametric rate refers to convergence with the rate $\varepsilon = n^{-1/2}$.) For further work on nonparametric estimation of change-points and discontinuous functions from direct observations, see, for example, [1, 9, 18, 21, 22, 25, 27] and the references cited therein. On the other hand, nonparametric estimation of a change-point from *indirect* observations, that is, for operators K that are not the identity, is much less studied. Furthermore, the literature contains some contradictory statements with regard to best achievable rates of convergence and therefore, leaves open the question of how to construct optimal estimators.

An important result in this area is due to Neumann [19], who investigated the problem of change-point estimation from indirect observations in a density deconvolution model. He assumes that the observations are $Y_i = X_i + \xi_i$, $i = 1, \ldots, n$, where $X_i$ are i.i.d. random variables with unknown probability density $f$ and where $\xi_i$ are i.i.d. random errors, independent of the $X_i$'s, with known probability density $K$. The problem considered by Neumann [19] is to estimate the location $\theta$ of a discontinuity jump in $f$, where this density is assumed to satisfy a Lipschitz condition away from the change-point. Neumann [19] proved that the order of the minimax risk in estimating $\theta$ is $\min\{n^{-2/(2\beta+3)}, n^{-1/(2\beta+1)}\}$, provided that the tails of the characteristic function $\widehat{K}(\omega)$ of $\xi_i$ decrease at the rate $|\omega|^{-\beta}$, $\beta > 0$. In the nonparametric regression context, Raimondo [21] considered the problem of estimating a change-point in the $\beta$th derivative of the regression function. Assuming that this derivative satisfies a Lipschitz condition apart from the change-point $\theta$, Raimondo [21] claims that the best rate of convergence in estimating $\theta$ is $n^{-1/(2\beta+1)}$. Estimation procedures that achieve this rate were also proposed by Wang [26] and, more recently, by Huh and Carrière [13] and Park and Kim [20]. Clearly, if $K$ is the Green's function of a linear differential operator of integer order $\beta$, estimating the change-point $\theta$ of $f$ from indirect observation as in model (1) is equivalent to estimating the change-point in the derivative of order $\beta$ from direct observations. This fact indicates that there is a discrepancy between the rates of convergence obtained by Neumann [19], on the one hand, and by Raimondo [21] and other authors cited above, on the other hand. In particular, the rates obtained by Neumann [19] are faster. Although asymptotic equivalence between the two indirect observation models (the density model as in [19], and regression/white noise model as in [21]) has not been established formally, it would seem natural to expect that the rates of convergence are in agreement. In what follows, we will show that the "faster" rates of Neumann [19] can indeed be attained and they are optimal for the white noise model (1). This fact will be deduced from more general results.

We study the problem of change-point estimation in model (1) for two different scales of functional classes $\mathcal{G}$ that quantify smoothness of $f$ away from the change-point. We derive lower bounds on the minimax risk (see Theorems 2 and 4) and develop rate optimal estimators (see Theorems 1 and 3).

In particular, we show that if $f$ can be represented as the sum of a jump function and a smooth function whose $m$th derivative exists and is bounded for all $x$, then the minimax risk in estimating $\theta$ is of order $\min\{\varepsilon^{(2m+2)/(2m+2\beta+1)}, \varepsilon^{2/(2\beta+1)}\}$, provided that the tails of the Fourier transform $\widehat{K}$ of $K$ behave like $|\omega|^{-\beta}$, as $|\omega| \to \infty$, with $\beta > 0$. The elbow in the rates of convergence corresponds to the cases where $\beta > 1/2$ and $0 < \beta \le 1/2$. If $\beta > 1/2$, the convolution kernel $K$ belongs to $\mathbb{L}_2(\mathbb{R})$. In what follows we call such convolution kernels and the corresponding setup *regular*. In contrast, under $0 < \beta \le 1/2$, the convolution kernel $K$ does not belong to $\mathbb{L}_2(\mathbb{R})$. We will call the latter case *singular* because it necessarily corresponds to a singular convolution integral in (2).

We introduce a new estimation technique that involves, as a key element, detection of zero crossings of an estimate of a properly smoothed *second derivative* of $f$. This differs from most change-point detection methods described in the statistical literature that typically use a properly smoothed first derivative of $f$. On the other hand, our second derivative based approach has parallels in digital image processing in the context of edge detection, where it is often referred to as the Laplacian method (cf. [11]). It is interesting to note that in the regular case seemingly intuitive procedures based on detecting a maximum in the first derivative lead to slower rates of convergence (see further discussion in Section 5).

The optimal rate of convergence in the regular case, $\varepsilon^{(2m+2)/(2m+2\beta+1)}$, clarifies how smoothness of $f$ away from the change-point (given by the index $m$) and ill-posedness of the kernel $K$ (given by $\beta$) affect achievable accuracy in change-point estimation from indirect observations. The result of Neumann [19] in the density deconvolution model, with standard calibration $\varepsilon = n^{-1/2}$, can be viewed as the "density analog" of a special case of our result with $m = 1$, that is, when $f$ is Lipschitz apart from the change-point. When the "smooth part" of the unknown function $f$ is analytic, our results show that in the regular case the optimal rate is $\varepsilon$, up to a logarithmic factor in $\varepsilon^{-1}$, that is, it is nearly the parametric rate. Interestingly, in this case the ill-posedness index $\beta$ of $K$ appears in the risk bound only as a power of the logarithmic factor. This means that ill-posedness of $K$ does not affect significantly the quality of estimation when $f$ is very smooth apart from the change-point. We also show that in the singular case the optimal rate of convergence is $\varepsilon^{2/(2\beta+1)}$, up to a logarithmic factor, regardless of the smoothness of $f$ away from the change-point.

Our results elucidate the following important feature of the problem: when estimating a change-point from indirect data, the best achievable rates of convergence depend on the behavior of the function $f$ away from the change-point location. This is in striking contrast to the direct observations case, where the rate $\varepsilon^2$ is the best one can achieve regardless of how many derivatives $f$ possesses apart from the discontinuity jump. Our results also indicate that the procedure of Raimondo [21] is not optimal when estimating a change-point of the $\beta$th derivative, $\beta \geq 1$, in the direct observation model, contrary to what is claimed in that paper (see further discussion in Section 5). We note that estimating a change-point from indirect observations can be done with higher accuracy than curve estimation in nonparametric deconvolution (see, e.g., [5, 6, 7, 10], and [8]).

The rest of the paper is organized as follows. Section 2 introduces notation and definitions of the functional classes. In Section 3 we construct a probe functional that is used for detection of the change-point from indirect observations; some properties of the probe functional are discussed, and its estimator is developed. Section 4 describes the two-stage change-point estimation procedure and presents our main results. Section 5 concludes with a discussion of the main results and Section 6 contains the proofs.

**2. Preliminaries.** We begin with some notation and definitions. Let $\widehat{g}$ or $(\mathsf{F}g)$ denote the Fourier transform of a function $g \in \mathbb{L}_2(\mathbb{R})$, in particular, if $g \in \mathbb{L}_1(\mathbb{R}) \cap \mathbb{L}_2(\mathbb{R})$,

$$\widehat{g}(\omega) \equiv (\mathsf{F}g)(\omega) \stackrel{\text{def}}{=} \int_{-\infty}^{\infty} g(x) e^{2\pi i \omega x} \, dx, \qquad \omega \in \mathbb{R}.$$

Let $f(x\pm) = \lim_{t \to x\pm} f(t)$ be the one-sided limits of $f$ at point $x$ and let $[f](x) = f(x+) - f(x-)$ be the local jump function. We say that $\theta \in \mathbb{R}$ is a *change-point* of $f$ if $[f](\theta) \neq 0$, and $f(\theta+)$ and $f(\theta-)$ are finite.

We will consider minimax estimation of a change-point $\theta$ of $f$ by assuming that $f$ belongs to one of the two functional classes, $\mathcal{F}_m$ or $\mathcal{A}_\nu$, defined below.

DEFINITION 1. Let $a, L > 0$ be fixed constants. We say that $f \in \mathcal{F}_1 = \mathcal{F}_1(a, L)$ if $f \in \mathbb{L}_2(\mathbb{R})$ and if $f$ has a single change-point $\theta \in [0, 1]$ such that $|[f](\theta)| \geq a$ and

$$|f(x) - f(x')| \leq L|x - x'| \qquad \forall x, x' \in \mathbb{R}, x \leq x', \theta \notin [x, x'].$$

The class $\mathcal{F}_1$ contains functions $f$ that have a single jump discontinuity of the first kind at $\theta \in [0, 1]$ and satisfy the Lipschitz condition on any interval that does not include $\theta$. This class was considered by Neumann [19] in the context of density deconvolution. To allow more smoothness of $f$ apart from the jump discontinuity, we introduce the following extension of $\mathcal{F}_1$.

DEFINITION 2. Let $a, L > 0$ and $m > 1$ be fixed constants. We say that $f \in \mathcal{F}_m = \mathcal{F}_m(a, L)$ if $f \in \mathbb{L}_2(\mathbb{R})$ and if $f$ has a single change-point $\theta \in [0, 1]$ such that the following conditions hold:

(i) We have $|[f](\theta)| \geq a$.

(ii) For all $x \neq \theta$ and $[f'](\theta) = 0$, $f'(x)$ exists so that the function $g_f : \mathbb{R} \to \mathbb{R}$ defined by

$$(3) \qquad g_f(x) = \begin{cases} f'(x), & x \neq \theta, \\ f'(\theta\pm), & x = \theta, \end{cases}$$

is continuous.

(iii) The function $g_f$ belongs to $\mathbb{L}_2(\mathbb{R})$ and its Fourier transform $\widehat{g}_f$ satisfies

$$(4) \qquad \int_{-\infty}^{\infty} |\widehat{g}_f(\omega)||\omega|^{m-1}\,d\omega \leq L.$$

If $m$ is an integer, condition (4) implies that the derivative $g_f^{(m-1)}$ exists and is bounded by $L$. In fact, (4) is only slightly stronger than this property. For example, (4) is valid if $g_f$ is in a Sobolev class of $\mathbb{L}_2$ smoothness $s > m - 1/2$, that is, when $\int |\widehat{g}_f(\omega)|^2 |\omega|^{2s}\,d\omega$ is bounded by an appropriate constant. This class is very close to, but smaller than, the class of functions with uniformly bounded derivative $g_f^{(m-1)}$.

It is important that in Definition 2 we have $[f'](\theta) = 0$. If $[f'](\theta) \neq 0$, parts (ii) and (iii) of Definition 2 cannot be satisfied, but we may still consider classes of functions $f$ that are smooth separately to the left and to the right of the change-point. However, introducing such classes seems to be unjustified, because additional smoothness in these terms does not improve the convergence rate of estimators of $\theta$. The minimax rate remains the same as for the class $\mathcal{F}_1$.

DEFINITION 3.   Let $\nu, a, L > 0$ be fixed constants. We say that $f \in \mathcal{A}_\nu = \mathcal{A}_\nu(a, L)$ if $f \in \mathbb{L}_2(\mathbb{R})$, and if $f$ has a single change-point $\theta \in [0, 1]$ such that conditions (i) and (ii) of Definition 2 are satisfied and

$$(5) \qquad \int_{-\infty}^{\infty} |\widehat{g}_f(\omega)|^2 \exp\{2\nu|\omega|\}\,d\omega \leq L^2,$$

where $g_f$ is defined in (3).

Assumption (5) implies that $g_f$ is infinitely differentiable and admits an analytical continuation onto a strip in the complex plane. Such classes of functions have been studied in the context of nonparametric estimation by many authors, starting with Ibragimov and Hasminskii [14]. For a recent overview, see, for example, [2].

The following assumption on $K$ will be used throughout this paper.

Assumption K. The function $K$ belongs to $\mathbb{L}_1(\mathbb{R})$, and there exist constants $\beta > 0$ (called the ill-posedness index of $K$) and $\underline{\kappa}, \overline{\kappa} > 0$, such that

$$(6) \qquad \underline{\kappa}(1 + |\omega|^2)^{-\beta/2} \leq |\widehat{K}(\omega)| \leq \overline{\kappa}(1 + |\omega|^2)^{-\beta/2} \qquad \forall \omega \in \mathbb{R}.$$

Assumption K is quite standard in deconvolution problems and corresponds to what is known as a *moderately ill-posed* problem. Green's functions of linear differential operators are important examples of kernels $K$ satisfying Assumption K for the regular case ($\beta > 1/2$). For instance, let $v(x) = e^x$, $-\infty < x < 0$, $v(0) = 1/2$ and $v(x) = 0$, $0 < x < \infty$. For nonvanishing real constants $b_j$, $j = 1, \ldots, k$, we define

$$(7) \qquad v_j(x) = |b_j| v(b_j x) \quad \text{and} \quad K = v_1 * v_2 * \cdots * v_k,$$

where $*$ stands for the convolution on $\mathbb{R}$. The Fourier transform of the kernel $K$ is given by $\widehat{K}(\omega) = \{\prod_{j=1}^k (1 - 2\pi b_j^{-1} i\omega)\}^{-1}$, and Assumption K holds with $\beta = k$. In this case $f$ can be recovered from $\mathsf{K}f$ by applying the linear differential operator

$$(8) \qquad f(x) = \left\{ \prod_{j=1}^k \left( 1 - 2\pi b_j^{-1} \frac{d}{dx} \right) \right\} (\mathsf{K}f)(x);$$

see [12], Chapter II. As for the singular case ($0 < \beta \leq 1/2$), examples are more peculiar; for instance, one may consider $K$ to be the probability density of a gamma distribution with shape parameter $\beta$.

## 3. Probe functional.

We will develop estimation procedures that are based on minimization of an empirical version of a properly chosen probe functional. Let $\varphi : \mathbb{R} \to \mathbb{R}$ be an even, twice continuously differentiable function that attains its global maximum at 0. Further conditions on $\varphi$ will be introduced below. Fix a bandwidth $h > 0$ and for $t \in \mathbb{R}$, $x \in \mathbb{R}$ define $\psi_t(x) = h^{-3} \varphi''(h^{-1}(x - t))$. Let $\langle \cdot, \cdot \rangle$ denote the standard inner product in $\mathbb{L}_2(\mathbb{R})$. Assuming that $\varphi'' \in \mathbb{L}_1(\mathbb{R}) \cap \mathbb{L}_2(\mathbb{R})$, we define the *probe functional*

$$(9) \qquad \begin{aligned} \ell_h(t) &= \langle f, \psi_t \rangle = \int_{-\infty}^{\infty} f(x) \psi_t(x) \, dx \\ &= \frac{1}{h^3} \int_{-\infty}^{\infty} f(x) \varphi'' \left( \frac{x - t}{h} \right) dx, \qquad t \in \mathbb{R}. \end{aligned}$$

The probe functional $\ell_h(t)$ is thus a smoothed second derivative of $f$ at point $t$: as $h$ tends to zero, $\ell_h(t)$ converges to $f''(t)$, provided that $f$ is twice continuously differentiable at $t$. The points $t$ where $|\ell_h(t)|$ is close to zero are indicative of the change-point location; this idea underlies the construction.

An estimator of $\ell_h(t)$ based on observations (1) can be developed as follows. Denote by $\mathsf{K}^*$ the adjoint operator to $\mathsf{K}$ given by

$$(10) \qquad (\mathsf{K}^*g)(x) = \int_{-\infty}^{\infty} K(y-x)g(y)\,dy, \qquad g \in \mathbb{L}_2(\mathbb{R}).$$

By the linear functional strategy, if $\psi_t \in \mathrm{Range}(\mathsf{K}^*)$, then there exists a function $\gamma_t \in \mathbb{L}_2(\mathbb{R})$ such that

$$\ell_h(t) = \langle f, \psi_t \rangle = \langle f, \mathsf{K}^*\gamma_t \rangle = \langle \mathsf{K}f, \gamma_t \rangle.$$

The function $\gamma_t$ satisfies $(\mathsf{K}^*\gamma_t)(x) = \psi_t(x) = h^{-3}\varphi''(h^{-1}(x-t))$ almost everywhere in $\mathbb{R}$. Taking the Fourier transforms, and using (10) and the fact that $\widehat{\psi_t}(\omega) = -(2\pi\omega)^2\widehat{\varphi}(\omega h)e^{2\pi i\omega t}$, we find

$$\widehat{\gamma_t}(\omega) = -(2\pi\omega)^2 e^{2\pi i\omega t} \frac{\widehat{\varphi}(\omega h)}{\overline{\widehat{K}(-\omega)}}.$$

We will always choose $\varphi$ so that $\widehat{\gamma_t} \in \mathbb{L}_1(\mathbb{R}) \cap \mathbb{L}_2(\mathbb{R})$. Under this assumption we may write

$$\gamma_t(x) = \int_{-\infty}^{\infty} \widehat{\gamma_t}(\omega)e^{-2\pi i\omega x}\,d\omega = -\int_{-\infty}^{\infty} (2\pi\omega)^2 e^{2\pi i\omega(t-x)} \frac{\widehat{\varphi}(\omega h)}{\overline{\widehat{K}(-\omega)}}\,d\omega$$

and

$$\ell_h(t) = \int_{-\infty}^{\infty} f(x)\psi_t(x)\,dx = \int_{-\infty}^{\infty} (\mathsf{K}f)(x)\gamma_t(x)\,dx.$$

Based on these considerations, we define the estimator $\tilde{\ell}_h(t)$ of $\ell_h(t)$ by

$$\tilde{\ell}_h(t) = \int_{-\infty}^{\infty} \gamma_t(x)\,dY(x), \qquad t \in \mathbb{R}.$$

Properties of estimation procedures that we develop are determined crucially by (i) accuracy of the probe functional estimator $\tilde{\ell}_h(t)$ and (ii) the ability of the probe functional to detect the change-point. The former is quantified by the next lemma, which we prove under the following assumption on the smoothing function $\varphi$.

ASSUMPTION 1.   The Fourier transform $\widehat{\varphi}$ of $\varphi \in \mathbb{L}_2(\mathbb{R})$ satisfies

$$\int_{-\infty}^{\infty} |\omega|^{2\beta+6}|\widehat{\varphi}(\omega)|^2\,d\omega < \infty;$$

that is, $\varphi$ belongs to the Sobolev space with smoothness index $\beta+3$.

LEMMA 1. *Let Assumption 1 and the left inequality in Assumption K hold, and let $h > 0$. Then the Gaussian random process $Z_h(t) = \tilde{\ell}_h(t) - \ell_h(t)$, $t \in B$, where $B$ is a subinterval of $[0, 1]$, satisfies*

$$(11) \qquad \mathbb{E}[Z_h(t)] \equiv 0, \qquad \sigma_Z^2 \stackrel{\text{def}}{=} \sup_{t \in B} \mathbb{E}[Z_h^2(t)] \leq C_1 \varepsilon^2 h^{-2\beta-5}.$$

*Furthermore, for any $\lambda \geq 2\sigma_Z$*

$$(12) \qquad \mathbb{P}\left\{ \sup_{t \in B} |Z_h(t)| \geq \lambda \right\} \leq C_2\left( \frac{h^{\beta+3/2}}{\varepsilon} \right)|B|\lambda \exp\left\{ -C_3 \frac{\lambda^2 h^{2\beta+5}}{\varepsilon^2} \right\},$$

*where $|B|$ stands for the Lebesgue measure of the set $B$ and $C_i$, $i = 1, 2, 3$, are positive constants.*

Further results will be obtained under the following condition which is stronger than Assumption 1.

ASSUMPTION 2. *The Fourier transform $\hat{\varphi}$ of $\varphi \in \mathbb{L}_2(\mathbb{R})$ is an even, nonnegative, infinitely differentiable function supported on $[-2/3, -1/3] \cup [1/3, 2/3]$ and taking the value 1 on $[-2/3 + \eta, -1/3 - \eta] \cup [1/3 + \eta, 2/3 - \eta]$ for some $\eta \in (0, 1/32)$.*

It follows from Assumption 2 that $\varphi$ is a real-valued, even, analytic function, rapidly decreasing at infinity, together with all its derivatives. In addition, since $\hat{\varphi}$ is nonnegative, $\varphi$ achieves its global maximum at $x = 0$, $\varphi'(0) = 0$ and $|\varphi''(0)| \geq M > 0$ for some constant $M$. The *Meyer wavelet* (see, e.g., [17], Section 7.2.2), centered at zero and rescaled accordingly, provides an example of a function that satisfies Assumption 2.

We summarize some properties of $\varphi'$ that will be repeatedly used in what follows.

(I) For all $x \in \mathbb{R}$, $\varphi'(0) = 0$ and $\varphi'(x) = -\varphi'(-x)$.

(II) The function $\varphi'$ decreases in $[0, 3/8]$, increases in $[3/4, 9/8]$ and has a unique minimum in $[3/8, 3/4]$, which is the point of the global minimum, $x = q_* \in [3/8, 3/4]$. By (I), $\varphi'$ attains its global maximum at $x = -q_* \in [-3/4, -3/8]$.

(III) There exists a unique zero of $\varphi'$ in the interval $[3/4, 3/2]$. We denote it by $q_0$ and let $d \stackrel{\text{def}}{=} |q_0 - q_*| > 0$. We have that $\varphi' < 0$ on $[0, q_0]$ and

$$(13) \qquad \inf_{x: |x - q_*| > d/2}\{\varphi'(x) - \varphi'(q_*)\} \geq r > 0$$

for some constant $r$.

Proofs of (I)–(III) are immediate and based on analysis of the integrand sign in the expressions for $\varphi'$ and $\varphi''$. Parameters $q_*$, $q_0$ and $d$ that appear in (II) and (III) depend on the specific choice of $\varphi$; condition (13) asserts that $q_*$ is a well-separated point of the global minimum of $\varphi'$.

The next lemma analyzes the separation between values of the probe functional $\ell_h(t)$ when $t$ varies in a "punctured" neighborhood of the change-point $\theta$.

LEMMA 2 (Separation rate). *Let Assumption* 2 *hold and let* $\delta \in (0, \bar{q}h)$, *where* $\bar{q} = q_* + 3d/4$, *and the constants* $q_*$ *and* $d$ *are given in* (II) *and* (III).

1. *Let* $f \in \mathcal{F}_m$ *and let*

$$(14) \qquad \delta \geq C_1(L/a)h^{m+1}$$

*for an absolute constant* $C_1 > 0$ *large enough. Then for sufficiently small* $h$,

$$(15) \qquad \inf_{t:\delta < |t-\theta| < \bar{q}h} \{|\ell_h(t)| - |\ell_h(\theta)|\} \geq C_2 a \delta h^{-3},$$

*where* $C_2$ *is a positive constant that depends on only* $a$, $L$ *and* $\varphi$.

2. *Let* $f \in \mathcal{A}_\nu$ *and let*

$$(16) \qquad \delta \geq C_3(L/a)h \exp\{-\nu/(3h)\}$$

*for an absolute constant* $C_3 > 0$ *large enough. Then* (15) *holds for sufficiently small* $h$.

The value that appears on the right-hand side (RHS) of (15) will be called the $\delta$-*separation rate* that corresponds to the probe functional $\ell_h$. Lemma 2 asserts that $\theta$ is a well-separated point of minimum of $|\ell_h(t)|$, provided that $h$ and $\delta$ satisfy (14) and (16) for $f \in \mathcal{F}_m$ and $f \in \mathcal{A}_\nu$, respectively. Conditions (14) and (16) are required to guarantee that the bias terms do not exceed the contrast expressed by the $\delta$-separation rate. They also show that if $f \in \mathcal{A}_\nu$, the value $\delta$ can be chosen much smaller than in the case of $f \in \mathcal{F}_m$; that is, the minimum of $|\ell_h(t)|$ is more pronounced when $f \in \mathcal{A}_\nu$. It is interesting to note that if a smoothed first derivative of $f$ is used as the probe functional and the maximum is sought, the corresponding $\delta$-separation rate is of order $\delta^2 h^{-3}$. As our proofs suggest, in the regular case this choice of the probe functional does not lead to a rate optimal estimation procedure (see Section 5).

**4. Estimation procedure and main results.** We are now in position to define the estimation procedure. The construction has two stages: First we localize the region that contains the change-point with probability close to 1 and then we search for a minimum of the absolute value of the probe functional inside the region.

The localization step is based on the following argument. As the proof of Lemma 2 shows, $\ell_h(t)$ is equal to $-h^{-2}\varphi'(h^{-1}(\theta - t))[f](\theta)$ up to a term that is negligible, provided that $|\theta - t| = O(h)$. It follows from (II) that $t_* \stackrel{\text{def}}{=} \arg\min_{t \in [0,1]} \{-\varphi'((\theta - t)h^{-1})[f](\theta)\}$ and $t^* \stackrel{\text{def}}{=} \arg\max_{t \in [0,1]} \{-\varphi'((\theta - t)h^{-1})[f](\theta)\}$ are within the distance $O(h)$ from $\theta$:

$$|t_* - \theta| = q_* h \in [3h/8, 3h/4], \qquad |t^* - \theta| = q_* h \in [3h/8, 3h/4].$$

In addition, $|t_* - t^*| = 2q_* h \in [3h/4, 3h/2]$. If $[f](\theta) < 0$, then $t^* > t_*$ and $[t_*, t^*]$ contains $\theta$; if $[f](\theta) > 0$, then $t^* < t_*$ and $\theta \in [t^*, t_*]$. Both $t_*$ and $t^*$ can be estimated from the data. This fact will be used to find an interval of size $O(h)$ that contains the change-point with probability close to 1.

Let

$$\hat{t}_* \stackrel{\text{def}}{=} \arg\min_{t \in [0,1]} \tilde{\ell}_h(t), \qquad \hat{t}^* \stackrel{\text{def}}{=} \arg\max_{t \in [0,1]} \tilde{\ell}_h(t) \tag{17}$$

and let $\hat{A}_h$ be the closed interval with endpoints $\hat{t}_*$ and $\hat{t}^*$. Our estimator $\tilde{\theta}_h$ of the change-point for given bandwidth $h$ is defined by

$$\tilde{\theta}_h = \arg\min_{t \in \hat{A}_h} |\tilde{\ell}_h(t)|. \tag{18}$$

We note that this construction depends on the bandwidth $h$ that will be chosen in an optimal way.

### 4.1. Functional class $\mathcal{F}_m$.

THEOREM 1. *Suppose that the left inequality in* (6) *is satisfied and that Assumption* 2 *holds. Let $\tilde{\theta}_*$ denote the change-point estimator $\tilde{\theta}_h$ with the bandwidth $h = h_*$ defined below.*

1. Regular case: *Assume that $\beta > 1/2$ and let*

$$h_* = C_1^* (\varepsilon/L)^{2/(2m+2\beta+1)}, \tag{19}$$

*where $C_1^* > 0$ is a constant. Then there exists a constant $C_2^* < \infty$ independent of $a, L$ such that*

$$R_\varepsilon[\tilde{\theta}_*; \mathcal{F}_m] \le C_2^* a^{-1} L^{(2\beta-1)/(2\beta+2m+1)} \varepsilon^{(2m+2)/(2m+2\beta+1)} \qquad \forall\, 0 < \varepsilon < 1.$$

2. Singular case: *Assume that $0 < \beta \le 1/2$ and let*

$$h_* = C_3^* \left(\varepsilon\sqrt{\ln\frac{1}{\varepsilon}}\right)^{2/(2\beta+1)}, \tag{20}$$

*where $C_3^* > 0$ is a sufficiently large constant. Then there exists a constant $C_4^* < \infty$ such that*

$$R_\varepsilon[\tilde{\theta}_*; \mathcal{F}_m] \le C_4^* \varepsilon^{2/(2\beta+1)} \left(\ln\frac{1}{\varepsilon}\right)^{(1-2\beta)/(2(2\beta+1))} \qquad \forall\, 0 < \varepsilon < 1.$$

THEOREM 2.   *Let the right-hand side inequality in* (6) *hold and assume that* $|\widehat{K}(\omega)| \neq 0$ *for all* $\omega$. *Then, for sufficiently small* $\varepsilon$, *the minimax risk over the class* $\mathcal{F}_m$ *is bounded from below as follows:*

1. Regular case: *If* $\beta > 1/2$, *then*

$$(21) \qquad R_\varepsilon^*[\mathcal{F}_m] \geq c_1^* a^{-1} L^{(2\beta-1)/(2\beta+2m+1)} \varepsilon^{(2m+2)/(2m+2\beta+1)},$$

*where* $c_1^*$ *does not depend on* $a, L$.

2. Singular case: *If* $0 \leq \beta \leq 1/2$, *then*

$$(22) \qquad R_\varepsilon^*[\mathcal{F}_m] \geq \begin{cases} c_2^* \varepsilon \left( \ln \dfrac{1}{\varepsilon} \right)^{-1/2}, & \text{if } \beta = 1/2, \\ c_3^* \varepsilon^{2/(2\beta+1)}, & \text{if } 0 < \beta < 1/2, \end{cases}$$

*where* $c_2^*, c_3^* > 0$ *are constants.*

Theorems 1 and 2 show that the estimate $\tilde{\theta}_*$ is rate optimal in the regular case. Moreover, the bounds identify the precise dependence of the minimax risk on the size of the jump $a$ and on the "Lipschitz constant" $L$ away from the jump (see definition of the class $\mathcal{F}_m$). As smoothness of $f$ away from the jump increases (i.e., as $m \to \infty$), the optimal rate of convergence approaches the usual parametric rate $\varepsilon$. In the singular case, $\tilde{\theta}_*$ is nearly rate optimal up to a factor logarithmic in $\varepsilon^{-1}$. Here the order of the rate of convergence is faster than the parametric rate $\varepsilon$, but slower than $\varepsilon^2$, the minimax rate achieved in change-point estimation under direct observations. It follows from Lemma 3 in Section 6 that the estimators $\hat{t}_*$ and $\hat{t}^*$ defined in (17) and associated with the bandwidth (20) are also nearly rate optimal in the singular case. We conjecture that for $0 < \beta < 1/2$ the extra logarithmic factor in the bound of Theorem 1 can be removed and thus $\varepsilon^{2/(2\beta+1)}$ is the optimal rate of convergence for such values of $\beta$.

### 4.2. Functional class $\mathcal{A}_\nu$.

THEOREM 3.   *Suppose that the left-hand side inequality in* (6) *is satisfied and that Assumption 2 holds. Let* $\tilde{\theta}_+$ *denote the change-point estimator* $\tilde{\theta}_h$ *with the bandwidth* $h = h_+$ *defined below.*

1. Regular case: *Assume that* $\beta > 1/2$ *and let*

$$(23) \qquad h_+ = \frac{\nu}{3} \left\{ \ln \frac{L}{\varepsilon} - \left( \beta + \frac{1}{2} \right) \ln \left( \frac{1}{\nu} \ln \frac{L}{\varepsilon} \right) \right\}^{-1}.$$

*Then there exists a constant* $C_5^* < \infty$ *independent of* $\nu, a, L$ *such that*

$$R_\varepsilon[\tilde{\theta}_+; \mathcal{A}_\nu] \leq C_5^* a^{-1} \varepsilon \left( \frac{1}{\nu} \ln \frac{L}{\varepsilon} \right)^{\beta-1/2} \qquad \forall\, 0 < \varepsilon < 1.$$

2. Singular case: *Assume that $0 < \beta \leq 1/2$ and let, for some sufficiently large $C_6^* > 0$,*

$$(24) \qquad h_+ = C_6^* \left( \varepsilon \sqrt{\ln \frac{1}{\varepsilon}} \right)^{2/(2\beta+1)}.$$

*Then there exists a constant $C_7^* < \infty$ such that*

$$R_\varepsilon[\tilde{\theta}_+; \mathcal{A}_\nu] \leq C_7^* \varepsilon^{2/(2\beta+1)} \left( \ln \frac{1}{\varepsilon} \right)^{(1-2\beta)/(2(2\beta+1))} \qquad \forall\, 0 < \varepsilon < 1.$$

THEOREM 4. *Let the right-hand side inequality in (6) hold. Then for sufficiently small $\varepsilon$ the minimax risk over the class $\mathcal{A}_\nu$ is bounded from below as follows:*

1. Regular case: *If $\beta > 1/2$, then*

$$(25) \qquad R_\varepsilon^*[\mathcal{A}_\nu] \geq c_4^* a^{-1} \varepsilon \left( \frac{1}{\nu} \ln \frac{L}{\varepsilon} \right)^{\beta - 1/2},$$

*where $c_4^*$ is a constant independent of $\nu, a, L$.*

2. Singular case: *If $0 \leq \beta \leq 1/2$, then*

$$(26) \qquad R_\varepsilon^*[\mathcal{A}_\nu] \geq \begin{cases} c_5^* \varepsilon \left( \ln \dfrac{1}{\varepsilon} \right)^{-1/2}, & \text{if } \beta = 1/2, \\ c_6^* \varepsilon^{2/(2\beta+1)}, & \text{if } 0 < \beta < 1/2, \end{cases}$$

*where $c_5^*, c_6^* > 0$ are constants independent of $\nu$.*

Theorems 3 and 4 indicate that $\tilde{\theta}_+$ is rate optimal in the regular case, and nearly rate optimal in the singular case. It is interesting to note that in the regular case, when $f \in \mathcal{A}_\nu$, almost parametric rates of convergence are attained by our estimation procedure. The ill-posedness of the convolution operator $\mathsf{K}$, as expressed by index $\beta$, does not have a significant effect on the rates of convergence when $f \in \mathcal{A}_\nu$; this fact is rather surprising.

**5. Discussion.** 1. Our technique elucidates how the construction of the probe functional affects estimation accuracy. An appropriate probe functional, $\ell_h$, should satisfy the following two requirements: (i) $\theta$ is a well-separated point of minimum (maximum) of $\ell_h(\cdot)$ or $|\ell_h(\cdot)|$; (ii) $\ell_h(t)$ admits a "good" estimator with "small" bias and variance. The proofs suggest that in the regular case the optimal rates of convergence are obtained by balancing three quantities: the $\delta$-separation rate, the bias and the stochastic error of estimation of a properly chosen probe functional. In the singular case the bias is asymptotically negligible and the optimal rates are obtained by balancing only two terms: the $\delta$-separation rate and the stochastic error.

As an illustration, consider, for instance, the regular case and estimation on classes $\mathcal{F}_m$. For our functional $\ell_h$, the $\delta$-separation rate in (15) is of order $\delta h^{-3}$, the stochastic error is characterized by the square root of the variance $\mathrm{Var}[\tilde{\ell}_h(t)] = O(\varepsilon^2 h^{-2\beta-5})$ [see (11)] and the bias term is of order $h^{m-2}$ [see (30) and (32)]. The balance between the three terms is given by the relationships $\delta h^{-3} \asymp h^{m-2} \asymp \varepsilon h^{-\beta-5/2}$. Solving for this, we get the optimal bandwidth $h \asymp \varepsilon^{2/(2m+2\beta+1)}$ and the corresponding optimal rate $\delta \asymp \varepsilon^{(2m+2)/(2m+2\beta+1)}$.

2. The proposed estimator is based on a local search for a zero of the smoothed second derivative of $f$. An alternative and seemingly natural approach would be to estimate $\theta$ by searching for a maximum of a smoothed first derivative of $f$, that is, to consider the probe functional $w_h(t) = h^{-2} \int f(x) \varphi'(h^{-1}(x-t)) \, dx$. This, however, does not lead to rate optimal estimation in the regular case. Although the stochastic error of the corresponding estimator $\tilde{w}_h(t)$ is smaller than that of $\tilde{\ell}_h(t)$ [$\mathrm{Var}[\tilde{w}_h(t)] = O(\varepsilon^2 h^{-2\beta-3})$], the $\delta$-separation rate is of order $\delta^2 h^{-3}$. The bias term is now of order $h^{m-1}$; this follows from similar arguments as in the proof of Lemma 2. By balancing the three terms (bias, stochastic error and $\delta$-separation rate), as explained in the previous remark in this section, it is not difficult to verify that the estimator based on a local maximization of $|\tilde{w}_h(t)|$ has risk of order $\varepsilon^{(m+2)/(2m+2\beta+1)}$ when $f \in \mathcal{F}_m$ and $\beta > 1/2$. We recall that the optimal rate of convergence given by Theorem 1 is faster, $\varepsilon^{(2m+2)/(2m+2\beta+1)}$.

3. The results of the present paper cover the problem of estimating change-points in the $\beta$th derivative of a function from direct observations. In particular, assume that $\beta$ is an integer and let $K$ be the Green's function of a linear differential operator of order $\beta$ as defined in (7). Denoting $q = \mathsf{K}f$, we note that estimating the change-point in $f$ from indirect observations (1) is equivalent to estimating the change-point in $q^{(\beta)}$ from *direct* observation of $q$ in the white noise model. Indeed, in view of the inversion formula (8), if $f$ has a change-point at $\theta$, then $q^{(\beta)}$ (or any linear differential form of $q$ of order $\beta$) has a change-point at $\theta$ as well. Therefore, regarding the observations of $q$ in the white noise model as indirect observations of $f$, with $K$ being the Green's function of a linear differential operator, we can apply our procedure to estimation of change-points in $q^{(\beta)}$. In particular, according to Theorems 1 and 2, if $q^{(\beta)}$ satisfies the Lipschitz condition away from the change-point (i.e., $m = 1$), the best achievable rate of convergence is $\varepsilon^{4/(2\beta+3)}$, which can be easily extended to the regression problem with equidistant design, where the rate becomes $n^{-2/(2\beta+3)}$. This indicates that change-point estimation procedures in [21], as well as in [13, 20, 26], are not rate optimal, contrary to what is claimed in some of these papers. Specifically, as a referee pointed out, the lower bound of Raimondo [21] is not correct. At the same time, our results are consistent with those obtained for density deconvolution by Neumann [19], who only considered the case of Lipschitz smoothness ($m = 1$).

**6. Proofs and auxiliary results.** In what follows $C$, $c$, $C_i$ and $c_i$, $i = 1, 2, \ldots$, stand for positive constants that may differ on different occurrences.

PROOF OF LEMMA 1. Assumptions 1 and K imply that $\gamma_t \in \mathbb{L}_2(\mathbb{R})$. Since $f \in \mathbb{L}_2(\mathbb{R})$ and $K \in \mathbb{L}_1(\mathbb{R})$ we have that $\mathsf{K}f \in \mathbb{L}_2(\mathbb{R})$; thus the Fourier transform $\widehat{\mathsf{K}f}$ exists and $\widehat{\mathsf{K}f} = \widehat{K}\widehat{f}$. Using (1) and Plancherel's formula we get, for any $t \in B$,

$$
\begin{aligned}
\mathbb{E}_f[\tilde{\ell}_h(t)] &= \int_{-\infty}^{\infty} \gamma_t(x)(\mathsf{K}f)(x)\,dx \\
&= \int_{-\infty}^{\infty} \widehat{\gamma}_t(\omega)\overline{\widehat{K}(\omega)\widehat{f}(\omega)}\,d\omega \\
&= -\int_{-\infty}^{\infty}(2\pi\omega)^2 e^{2\pi i\omega t}\widehat{\varphi}(\omega h)\overline{\widehat{f}(\omega)}\,d\omega \\
&= \int_{-\infty}^{\infty} \widehat{\psi}_t(\omega)\overline{\widehat{f}(\omega)}\,d\omega = \langle f, \psi_t \rangle = \ell_h(t),
\end{aligned}
$$

which proves that $\mathbb{E}[Z_h(t)] = 0$. Thus $Z_h(t)$ is a zero mean Gaussian random variable with variance

$$
\begin{aligned}
\mathbb{E}[Z_h^2(t)] &= \mathbb{E}\left[\left(\varepsilon\int_{-\infty}^{\infty}\gamma_t(x)\,dW(x)\right)^2\right] \\
&= \varepsilon^2\int_{-\infty}^{\infty}|\gamma_t(x)|^2\,dx \\
&= \varepsilon^2\int_{-\infty}^{\infty}|2\pi\omega|^4\frac{|\widehat{\varphi}(\omega h)|^2}{|\widehat{K}(-\omega)|^2}\,d\omega \\
&\leq c_1\frac{\varepsilon^2}{h^{2\beta+5}}\int_{-\infty}^{\infty}|\widehat{\varphi}(\omega)|^2(1+|\omega|^2)^{\beta+2}\,d\omega \\
&\leq c_2\varepsilon^2 h^{-2\beta-5},
\end{aligned}
$$

where we have used Assumptions K and 1. This proves (11).

To prove (12) we apply the inequality on the tails of Gaussian processes (see, e.g., [24], Proposition A.2.7). For this purpose, using Plancherel's formula and Assumptions K and 1 we obtain, for $t, s \in [0, 1]$,

$$
\begin{aligned}
\sigma^2(t, s) &\stackrel{\text{def}}{=} \mathbb{E}|Z(t) - Z(s)|^2 \\
&= \varepsilon^2\int_{-\infty}^{\infty}|\gamma_t(x) - \gamma_s(x)|^2\,dx \\
&= \varepsilon^2\int_{-\infty}^{\infty}|2\pi\omega|^4\frac{|\widehat{\varphi}(\omega h)|^2}{|\widehat{K}(-\omega)|^2}|e^{2\pi i\omega t} - e^{2\pi i\omega s}|^2\,d\omega
\end{aligned}
$$

$$\leq c_3 \varepsilon^2 |t-s|^2 \int_{-\infty}^{\infty} |\omega|^6 (1+|\omega|^2)^\beta |\widehat{\varphi}(\omega h)|^2 \, d\omega$$

$$\leq c_3 \varepsilon^2 h^{-2\beta-7} |t-s|^2 \int_{-\infty}^{\infty} (1+|\omega|^2)^{\beta+3} |\widehat{\varphi}(\omega)|^2 \, d\omega$$

$$\leq c_4 \varepsilon^2 h^{-2\beta-7} |t-s|^2.$$

Therefore, the number of balls of radius $r$ in the seminorm $\sigma(t, s)$ that cover the interval $B \subseteq [0, 1]$ does not exceed $c_5 r^{-1} \varepsilon h^{-\beta-7/2} |B|$, and applying Proposition A.2.7 of [24] (putting in the notation of that proposition $\varepsilon_0 = \sigma_Z$, $K = \varepsilon h^{-\beta-7/2} |B|$), we get the lemma.  $\square$

PROOF OF LEMMA 2.  1. Fix $t$ satisfying $\delta < |t - \theta| < \bar{q}h$ and define $\tau = (\theta - t)/h$; clearly $\bar{q} > |\tau| > \delta/h$. By (9),

$$
\begin{aligned}
(27) \quad \ell_h(t) &= \frac{1}{h^3} \int_{-\infty}^{\theta} \varphi'' \left( \frac{x-t}{h} \right) f(x) \, dx + \frac{1}{h^3} \int_{\theta}^{\infty} \varphi'' \left( \frac{x-t}{h} \right) f(x) \, dx \\
&= \frac{1}{h^2} \int_{-\infty}^{\tau} \varphi''(x) f(t+xh) \, dx + \frac{1}{h^2} \int_{\tau}^{\infty} \varphi''(x) f(t+xh) \, dx.
\end{aligned}
$$

First assume that $m = 1$. Let

$$
\begin{aligned}
J_1(t) &\stackrel{\text{def}}{=} \frac{1}{h^2} \int_{-\infty}^{\tau} \varphi''(x) [f(t+xh) - f(\theta-)] \, dx \\
&\quad + \frac{1}{h^2} \int_{\tau}^{\infty} \varphi''(x) [f(t+xh) - f(\theta+)] \, dx.
\end{aligned}
$$

Then using Definition 1 and the fact that $\varphi'(-\infty) = \varphi'(\infty) = 0$, we obtain

$$
\begin{aligned}
(28) \quad \ell_h(t) &= \frac{1}{h^2} \varphi'(\tau) f(\theta-) - \frac{1}{h^2} \varphi'(\tau) f(\theta+) + J_1(t) \\
&= -\frac{1}{h^2} \varphi'(\tau) [f](\theta) + J_1(t).
\end{aligned}
$$

Recall that $\varphi'(0) = 0$ and $|\varphi''(0)| \geq M > 0$. In addition, by (I)–(III), $\varphi'(x)$ has a unique zero in the interval $[-\bar{q}, \bar{q}]$ at $x = 0$. Therefore, $|\varphi'(\tau)| > c_1 |\tau|$ for all $|\tau| \in (\delta/h, \bar{q})$ and, for $h$ small enough, we get

$$(29) \quad \frac{1}{h^2} |[f](\theta)| |\varphi'(\tau)| \geq c_2 a h^{-2} |\tau| \geq c_3 a \delta h^{-3} .$$

Furthermore, we note that $\ell_h(\theta) = J_1(\theta)$ and that, for all $t$ satisfying $\delta < |t - \theta| < \bar{q}h$,

$$(30) \quad |J_1(t)| \leq \frac{L}{h} \int_{-\infty}^{\infty} |\varphi''(x)| |\tau - x| \, dx \leq c_4 L h^{-1}.$$

Using (14) we conclude that for sufficiently small $h$, the sign of $\ell_h(t)$ is determined by the first term in (28). Therefore, (15) holds for $m = 1$.

If $m > 1$, then integrating by parts in both integrals on the RHS of (27) and using the fact that $\varphi'(-\infty) = \varphi'(\infty) = 0$, we obtain

$$(31) \qquad \ell_h(t) = -\frac{1}{h^2}\varphi'(\tau)[f](\theta) - J_2(t)$$

with

$$J_2(t) \overset{\text{def}}{=} \frac{1}{h}\int_{-\infty}^{\infty}\varphi'(x)g_f(t+xh)\,dx$$
$$= \frac{1}{h^2}\int_{-\infty}^{\infty}(-2\pi i\omega)\widehat{\varphi}(\omega)e^{2\pi i\omega t/h}\widehat{g}_f(-\omega/h)\,d\omega,$$

where $g_f$ is defined in (3) and the last equality follows from the Plancherel formula. In view of Assumption 2 and Definition 2,

$$(32) \quad |J_2(t)| \le \frac{2\pi}{h^2}\max_{1/3 \le |\omega| \le 2/3}\frac{|\widehat{\varphi}(\omega)|}{|\omega|^{m-2}}\int_{-\infty}^{\infty}|\omega|^{m-1}|\widehat{g}_f(\omega/h)|\,d\omega \le c_5 L h^{m-2}.$$

This along with (14), (29), (31) and the fact that $\ell_h(\theta) = -J_2(\theta)$ completes the proof of the first statement of the lemma.

2. If $f \in \mathcal{A}_\nu$, then $\ell_h(t)$ is again given by (31) and $J_2(t)$ is now bounded using the Cauchy–Schwarz inequality:

$$|J_2(t)| \le 4\pi^2\left\{\int_{-\infty}^{\infty}|\omega|^2|\widehat{\varphi}(\omega h)|^2\exp\{-2\nu|\omega|\}\,d\omega\right\}^{1/2}$$
$$\times \left\{\int_{-\infty}^{\infty}|\widehat{g}_f(\omega)|^2\exp\{2\nu|\omega|\}\,d\omega\right\}^{1/2}$$
$$\le c_6\left\{\int_{1/(3h) \le |\omega| \le 2/(3h)}|\omega|^2\exp\{-2\nu|\omega|\}\,d\omega\right\}^{1/2}$$
$$\le c_7 L h^{-2}\exp\{-\nu/(3h)\}.$$

The same considerations as above complete the proof. $\square$

LEMMA 3. *Let Assumption 1 and the left inequality in Assumption K hold, and let $h^{\beta+1/2}\varepsilon^{-1} \ge C_1$. Then, for $f \in \mathcal{F}_m$ or $f \in \mathcal{A}_\nu$ and for all $h$ small enough,*

$$\max\{\mathbb{P}_f\{|\hat{t}_* - t_*| > hd/2\}, \mathbb{P}_f\{|\hat{t}^* - t^*| > hd/2\}\}$$
$$\le C_2 h^{\beta-1/2}\varepsilon^{-1}\exp\left\{-C_3\frac{h^{2\beta+1}}{\varepsilon^2}\right\},$$

*where $d > 0$ is given in* (III).

PROOF. We will derive the inequality for $\mathbb{P}_f\{|\hat{t}_* - t_*| > hd/2\}$ only; the proof of the other part is identical in every detail. Define $\Delta = \{t \in [0,1] : |t - t_*| > hd/2\}$. By definition of $t_*$ and $\hat{t}_*$, we have

$$\mathbb{P}_f\{|\hat{t}_* - t_*| > hd/2\} \leq \mathbb{P}_f\{\exists t \in \Delta : \tilde{\ell}_h(t_*) \geq \tilde{\ell}_h(t)\}$$
$$= \mathbb{P}_f\{\exists t \in \Delta : [\tilde{\ell}_h(t_*) - \ell_h(t_*)] + [\ell_h(t) - \tilde{\ell}_h(t)]$$
$$\geq \ell_h(t) - \ell_h(t_*)\}$$
$$\leq \mathbb{P}_f\left\{2 \sup_{t \in [0,1]} |\tilde{\ell}_h(t) - \ell_h(t)| \geq \inf_{t \in \Delta}(\ell_h(t) - \ell_h(t_*))\right\}.$$

It follows from (28) and (31) that

$$\ell_h(t) - \ell_h(t_*) = -\frac{1}{h^2}[f](\theta)\left\{\varphi'\left(\frac{\theta - t}{h}\right) - \varphi'\left(\frac{\theta - t_*}{h}\right)\right\} + J(t) - J(t_*),$$

where $J(t) \stackrel{\text{def}}{=} J_1(t)$ if $f \in \mathcal{F}_1$ and $J(t) \stackrel{\text{def}}{=} -J_2(t)$ if $f \in \mathcal{F}_m$, $m > 1$ or $f \in \mathcal{A}_\nu$, with $J_1(t)$ and $J_2(t)$ as defined in the proof of Lemma 2. Therefore, we obtain

$$\inf_{t \in \Delta}(\ell_h(t) - \ell_h(t_*))$$

$$(33) \qquad \geq \inf_{t \in \Delta} \frac{1}{h^2}[f](\theta)\left\{\varphi'\left(\frac{\theta - t_*}{h}\right) - \varphi'\left(\frac{\theta - t}{h}\right)\right\} - 2\sup_{t \in \mathbb{R}} |J(t)|$$

$$\geq \frac{c_1}{h^2} \inf_{x : |x - q_*| > d/2}\{\varphi'(x) - \varphi'(q_*)\} - 2\sup_{t \in \mathbb{R}} |J(t)|,$$

where the last inequality follows by change of variables, by definition of $t_*$ and because $\varphi'(x) = -\varphi'(-x)$. Using property (III), we obtain that the first term on the RHS of (33) is at least $c_1 r h^{-2}$, while the second one does not exceed in absolute value $c_2 L h^{m-2}$ if $f \in \mathcal{F}_m$ and $c_3 L h^{-2} \exp\{-\nu/(3h)\}$ if $f \in \mathcal{A}_\nu$ (see the proof of Lemma 2, where upper bounds on $|J(t)|$ were established). Noting that $h^{-2}$ dominates both $h^{m-2}$ and $h^{-2} \exp\{-\nu/(3h)\}$ as $h$ tends to zero and applying Lemma 1, we obtain that

$$\mathbb{P}_f\{|\hat{t}_* - t_*| > hd/2\} \leq \mathbb{P}_f\left\{\sup_{t \in [0,1]} |\tilde{\ell}_h(t) - \ell_h(t)| \geq c_4 h^{-2}\right\}$$

$$\leq c_5 h^{\beta - 1/2} \varepsilon^{-1} \exp\left\{-c_6 \frac{h^{2\beta+1}}{\varepsilon^2}\right\},$$

as claimed. □

PROOF OF THEOREM 1. 1. We begin with the regular case. The choice of $h = h_*$ in (19) implies that $\varepsilon^{-1} h^{\beta+1/2} \asymp \varepsilon^{-m/(m+\beta+1/2)}$, so Lemma 3 can be applied. Let $\Omega$ be the event that $|\hat{t}_* - t_*| \leq hd/2$ and $|\hat{t}^* - t^*| \leq hd/2$.

Recall that $|t_* - \theta| = |t^* - \theta| = q_* h$, where $q_*$ is defined in (II). Therefore, on the set $\Omega$,

$$(34) \qquad |\hat{t}_* - \theta| \le q_* h + hd/2 \le \bar{q}h \quad \text{and} \quad |\hat{t}^* - \theta| \le q_* h + hd/2 \le \bar{q}h,$$

where $\bar{q}$ is defined in Lemma 2. Recall that, by property (III) $\varphi'(x)$ has a unique zero in the interval $[-\bar{q}, \bar{q}]$ at the point $x = 0$. This guarantees that if $\Omega$ holds, the set $\hat{A}_h$ contains a unique zero of the function $t \mapsto \varphi'(h^{-1}(\theta - t))$ at $t = \theta$ and thus the definition (18) is justified.

We write

$$(35) \qquad \begin{aligned} \mathbb{E}_f |\tilde{\theta}_h - \theta|^2 &= \mathbb{E}_f\{|\tilde{\theta}_h - \theta|^2 \mathbf{1}(\Omega)\} + \mathbb{E}_f\{|\tilde{\theta}_h - \theta|^2 \mathbf{1}(\Omega^c)\} \\ &\le \mathbb{E}_f\{|\tilde{\theta}_h - \theta|^2 \mathbf{1}(\Omega)\} + \mathbb{P}_f(\Omega^c), \end{aligned}$$

where $\mathbf{1}(\cdot)$ denotes the indicator function. By Lemma 3 we have

$$(36) \qquad \begin{aligned} \mathbb{P}(\Omega^c) &= \mathbb{P}\{(|\hat{t}_* - t_*| > hd/2) \cup (|\hat{t}^* - t^*| > hd/2)\} \\ &\le c_1 h^{\beta+1/2} \varepsilon^{-1} \exp\{-c_2 h^{2\beta+1} \varepsilon^{-2}\} \\ &\le c_3 \varepsilon^{-m/(m+\beta+1/2)} \exp\{-c_4 \varepsilon^{-(2m)/(m+\beta+1/2)}\}. \end{aligned}$$

Furthermore, when $\Omega$ holds, it follows from (34) and from the construction of $\tilde{\theta}_h$ that $|\tilde{\theta}_h - \theta| \le \bar{q}h$. Thus for any $\delta \in (0, \bar{q}h)$, the first term on the RHS of (35) can be bounded as

$$(37) \qquad \mathbb{E}_f\{|\tilde{\theta}_h - \theta|^2 \mathbf{1}(\Omega)\} \le \delta^2 + \sum_{j=1}^{J} \delta^2 2^{2j} \mathbb{P}_f(\{|\tilde{\theta}_h - \theta| \in \Delta_j\} \cap \Omega),$$

where $\Delta_j \stackrel{\text{def}}{=} [\delta 2^{j-1}, \delta 2^j]$ and $J = \min\{j : \delta 2^j > \bar{q}h\}$. Let $T_j = \{t : |t - \theta| \in \Delta_j\}$; we note that $|T_j| = \delta 2^{j-1}$. Then we have

$$(38) \qquad \begin{aligned} &\mathbb{P}_f(\{|\tilde{\theta}_h - \theta| \in \Delta_j\} \cap \Omega) \\ &\le \mathbb{P}_f(\exists t \in T_j : |\tilde{\ell}_h(\theta)| \ge |\tilde{\ell}_h(t)|) \\ &\le \mathbb{P}_f\left\{2 \sup_{t \in T_j} |\tilde{\ell}_h(t) - \ell_h(t)| \ge \inf_{t \in T_j}(|\ell_h(t)| - |\ell_h(\theta)|)\right\}. \end{aligned}$$

We first estimate $\inf_{t \in T_j}(|\ell_h(t)| - |\ell_h(\theta)|)$ using Lemma 2. Note that Lemma 2 can be applied with $\{t : \delta \le |t - \theta| \le \bar{q}h\}$ replaced by $T_j$ for each $j = 1, \ldots, J$, provided that $\delta 2^{j-1} \ge c_5(L/a)h^{m+1}$. In particular, we have

$$(39) \qquad \begin{aligned} \inf_{t \in T_j}(|\ell_h(t)| - |\ell_h(\theta)|) &\ge \inf_{t \,:\, \delta 2^{j-1} < |t - \theta| < \bar{q}h}(|\ell_h(t)| - |\ell_h(\theta)|) \\ &\ge c_6 a\delta h^{-3} 2^j, \qquad j = 1, \ldots, J. \end{aligned}$$

Let

$$(40) \quad \delta = c_7 (L/a) h^{m+1} = c_8 L^{(2\beta-1)/(2\beta+2m+1)} a^{-1} (\varepsilon^2)^{(m+1)/(2\beta+2m+1)}.$$

It is straightforward to verify that with this choice of $\delta$ and $h$ for sufficiently large $c_7$, conditions of Lemma 2 are satisfied. In addition, $2^j a \delta h^{-3} \geq c_9 \varepsilon h^{-\beta-5/2}$ for some constant $c_9$ and each $j = 1, \ldots J$. Therefore, using (15) and applying Lemma 1, we obtain

$$\mathbb{P}_f \left\{ \sup_{t \in T_j} |Z_h(t)| \geq c_6 a \delta h^{-3} 2^j \right\}$$

$$(41) \qquad \leq c_{10} \left( \frac{h^{\beta+3/2}}{\varepsilon} \right) |T_j| a \delta h^{-3} 2^j \exp \left\{ -c_{11} \frac{(a \delta h^{-3} 2^j)^2 h^{2\beta+5}}{\varepsilon^2} \right\}$$

$$\qquad \leq c_{12} \frac{h^{\beta+3/2}}{\varepsilon} a \delta^2 2^{2j} h^{-3} \exp \left\{ -c_{13} 2^{2j} \frac{a^2 \delta^2 h^{2\beta-1}}{\varepsilon^2} \right\},$$

$$\qquad \leq c_{12} \frac{h^{\beta-3/2}}{\varepsilon} a \delta^2 2^{2j} \exp \{-c_{14} 2^{2j}\},$$

where we have taken into account that $\delta^2 h^{2\beta-1} \varepsilon^{-2} \geq c > 0$ under (19) and (40). Note also that $c_{14}$ can be made large enough by choice of $c_8$ in (40). Furthermore, since $h^{\beta-3/2} \delta^2 \varepsilon^{-1} = h^{\beta+2m+1/2} \varepsilon^{-1} = \varepsilon^{m/(\beta+m+1/2)} = o(1)$ as $\varepsilon \to 0$, we finally obtain from (41), (38) and (37) that

$$\mathbb{E}_f \{ |\tilde{\theta}_h - \theta|^2 \mathbf{1}(\Omega) \} \leq \delta^2 + \delta^2 o(1) \sum_{j=1}^J 2^{2j} \exp\{-c_{14} 2^{2j}\}$$

$$\leq \delta^2 (1 + o(1)).$$

Combining this with (35) and (36), we complete the proof of the first part of the theorem.

2. For the singular case, the proof follows the same lines with minor modifications; we indicate them below.

The choice of $h$ in (20) ensures that $h^{\beta+1/2} \varepsilon^{-1} = c_{15} \sqrt{\ln \varepsilon^{-1}}$ so that Lemma 3 applies. In addition, by choice of $C_3^*$ large enough, the probability $\mathbb{P}(\Omega^c) = o(h^{2+\eta})$ for any $\eta > 0$ and $\varepsilon \to 0$. Arguing as in the proof of the first part, we see that inequalities (37)–(39) hold. For some constant $c_{16} > 0$, let

$$(42) \qquad \delta = c_{16} \varepsilon h^{-\beta+1/2} = c_{17} \varepsilon^{1/(\beta+1/2)} \left( \sqrt{\ln \frac{1}{\varepsilon}} \right)^{(-\beta+1/2)/(\beta+1/2)}.$$

Under this choice,

$$\inf_{t \in T_j} (|\ell_h(t)| - |\ell_h(\theta)|) \geq c_{18} \delta h^{-3} 2^j$$

$$\geq c_{19} \varepsilon h^{-\beta-5/2} 2^j, \qquad j = 1, \ldots, J.$$

The last inequality ensures that Lemma 1 can be applied and, similarly to (41), we have

$$\mathbb{P}_f\left\{\sup_{t\in T_j}|Z_h(t)|\geq c_{18}\delta h^{-3}2^j\right\}\leq c_{20}\frac{h^{\beta-3/2}}{\varepsilon}\delta^2 2^{2j}\exp\left\{-c_{21}\frac{\delta^2 2^{2j}h^{2\beta-1}}{\varepsilon^2}\right\}$$

$$\leq c_{22}\frac{\varepsilon}{h^{\beta+1/2}}2^{2j}\exp\{-c_{23}2^{2j}\}.$$

Substituting expression (20) for $h$ and summing up over $j = 1, \ldots, J$, we complete the proof. $\square$

PROOF OF THEOREM 2. We use the method of proving minimax lower bounds based on a reduction to the problem of testing two simple hypotheses; see, for example, [16], Chapter 2, or [23], Chapter 2.

Pick $f_0 \in \mathcal{F}_m(a, L/2)$ such that $f_0$ has a unique jump discontinuity at $\theta_0 = 0$ and $[f_0](0) = a$. Fix $\delta \in (0, 1]$ and define $v(x) = a\mathbf{1}_{[0,\delta]}(x)$, $x \in \mathbb{R}$. The Fourier transform of $v$ is given by

$$\widehat{v}(\omega) = a\int_0^\delta e^{2\pi i\omega x}\,dx = \frac{a}{2\pi i\omega}[e^{2\pi i\omega\delta}-1].$$

Fix $N > 0$ and define

$$v_N(x) = \int_{-N}^N \widehat{v}(\omega)e^{-2\pi ix\omega}\,d\omega, \qquad x \in \mathbb{R}.$$

The Fourier transform of this function is $\widehat{v}_N(\omega) = \widehat{v}(\omega)\mathbf{1}(|\omega| \leq N)$. Let

$$f_1(x) = f_0(x) - [v(x) - v_N(x)], \qquad x \in \mathbb{R}.$$

The function $x \mapsto f_0(x) - v(x)$ has a unique jump at $x = \delta$ with $[f_0 - v](\delta) = -a$, while $v_N$ is infinitely differentiable. Therefore, $x \mapsto f_1(x)$ has a unique jump at $x = \delta$ with $[f_1](\delta) = -a$. Set $\theta_1 = \delta$, where the index 1 indicates that $\theta_1$ is the change-point of $f_1$.

We now show that $f_1 \in \mathcal{F}_m(a, L)$ under appropriate choice of $N$. First, clearly $f_1 \in \mathbb{L}_2(\mathbb{R})$, since $f_0, v, v_N \in \mathbb{L}_2(\mathbb{R})$. Next, the Fourier transform of the derivative $v_N'(x)$ is given by $(\mathsf{F}v_N')(\omega) = (-2\pi i\omega)\widehat{v}_N(\omega) = (-2\pi i\omega)\widehat{v}(\omega)\mathbf{1}(|\omega| \leq N)$ and

$$\begin{aligned}
\int_{-\infty}^\infty |(\mathsf{F}v_N')(\omega)||\omega|^{m-1}\,d\omega &= 2\pi\int_{-\infty}^\infty |\widehat{v}_N(\omega)||\omega|^m\,d\omega \\
&= a\int_{-N}^N |e^{2\pi i\omega\delta}-1||\omega|^{m-1}\,d\omega \\
&\leq \frac{4\pi a}{m+1}\delta N^{m+1}.
\end{aligned}$$

(43)

In what follows choose $N = (\frac{(m+1)L}{8\pi a\delta})^{1/(m+1)}$. Then the expression in (43) is less than $L/2$.

First let $m = 1$. Then (43) implies that the derivative $|v'_N(x)|$ is uniformly in $x \in \mathbb{R}$ bounded by $L/2$ and thus $v_N$ is Lipschitz continuous with Lipschitz constant $L/2$ on $\mathbb{R}$. Also, $f_0 - v$ has this property apart from $\theta_1 = \delta$. Hence, $f_1$ is Lipschitz continuous with Lipschitz constant $L$ apart from $\theta_1 = \delta$, which proves that $f_1 \in \mathcal{F}_1(a, L)$.

Now let $m > 1$. Then $g_{f_1} = g_{f_0} + v'_N$, with $g_f$ defined in (3) and

$$\int_{-\infty}^{\infty} |\hat{g}_{f_0}(\omega)||\omega|^{m-1} \, d\omega \le L/2,$$

since $f_0 \in \mathcal{F}_m(a, L/2)$. This and (43) prove that, under our choice of $N$,

$$\int_{-\infty}^{\infty} |\hat{g}_{f_1}(\omega)||\omega|^{m-1} \, d\omega \le L.$$

We have thus shown that $f_1 \in \mathcal{F}_m(a, L)$.

For brevity, let $\mathbb{P}_0$ and $\mathbb{P}_1$ denote the probability measures associated with the observations $\mathbf{Y} = \{Y(x) : x \in \mathbb{R}\}$ in model (1) with $f = f_0$ and $f = f_1$, respectively. In view of Girsanov's formula, the Kullback–Leibler divergence between $\mathbb{P}_0$ and $\mathbb{P}_1$ has the form

$$\mathcal{K}(\mathbb{P}_0, \mathbb{P}_1) \stackrel{\text{def}}{=} \int \ln \frac{d\mathbb{P}_0}{d\mathbb{P}_1} \, d\mathbb{P}_0 = \frac{1}{2\varepsilon^2} \|\mathsf{K}(f_0 - f_1)\|_2^2.$$

The function $\Delta = f_0 - f_1 = v - v_N$ belongs to $\mathbb{L}_2(\mathbb{R})$ and its Fourier transform is given by $\hat{v}_N(\omega)\mathbf{1}(|\omega| > N)$. Since $K \in \mathbb{L}_1(\mathbb{R})$, $\widehat{\mathsf{K}\Delta}$ exists and $\widehat{\mathsf{K}\Delta} = \widehat{K}\hat{\Delta}$. Hence, by Plancherel's formula,

$$(44) \qquad \mathcal{K}(\mathbb{P}_0, \mathbb{P}_1) = \frac{1}{2\varepsilon^2} \int_{|\omega| > N} |\hat{v}(\omega)|^2 |\widehat{K}(\omega)|^2 \, d\omega.$$

Assume first that $\beta > 1/2$. Then

$$\mathcal{K}(\mathbb{P}_0, \mathbb{P}_1) \le c_1 \frac{(a\delta)^2}{\varepsilon^2} N^{-2\beta+1}$$

$$= c_2 \varepsilon^{-2} L^{(1-2\beta)/(m+1)} (a\delta)^{(2\beta+2m+1)/(m+1)},$$

where we have used Assumption K and the fact that $|\hat{v}(\omega)| \le a\delta \; \forall \omega \in \mathbb{R}$. Choosing

$$\delta \asymp a^{-1} L^{(2\beta-1)/(2\beta+2m+1)} \varepsilon^{2(m+1)/(2\beta+2m+1)},$$

we ensure that $\mathcal{K}(\mathbb{P}_0, \mathbb{P}_1) \le \alpha < \infty$ for $\varepsilon$ small enough. On the other hand, $|\theta_0 - \theta_1| = \delta$ and it follows from part (iii) of Theorem 2.2 in [23] that $\sup_{f \in \mathcal{F}_m} \mathbb{E}_f |\hat{\theta} - \theta|^2 \ge c_3 \delta^2$. This completes the proof of (21) for $\beta > 1/2$.

Now let $0 < \beta \leq 1/2$. We decompose the domain of integration in (44) into two parts: $N < |\omega| \leq N'$ and $|\omega| > N'$, where $N' = 1/\delta$. For $N < |\omega| \leq N'$ we bound the integrand as above, while for $|\omega| > N'$ we use the fact that $|\hat{v}(\omega)| \leq (\pi|\omega|)^{-1}$. This yields

$$(45) \qquad \mathcal{K}(\mathbb{P}_0, \mathbb{P}_1) \leq c_4 \left( \frac{\delta^2}{\varepsilon^2} \int_{N < |\omega| \leq N'} \frac{d\omega}{|\omega|^{2\beta}} + \frac{1}{\varepsilon^2} \int_{|\omega| > N'} \frac{d\omega}{|\omega|^{2+2\beta}} \right).$$

Hence, for $\beta = 1/2$ we obtain

$$\mathcal{K}(\mathbb{P}_0, \mathbb{P}_1) \leq c_5 \left( \frac{\delta^2}{\varepsilon^2} \ln N' + \frac{1}{(\varepsilon N')^2} \right),$$

and the choice of $\delta \asymp \varepsilon (\ln \frac{1}{\varepsilon})^{-1/2}$ allows us to conclude the proof using the same argument as in the case of $\beta > 1/2$. Finally, for $0 < \beta < 1/2$, we get from (45) that

$$\mathcal{K}(\mathbb{P}_0, \mathbb{P}_1) \leq c_6 \left( \frac{\delta^2 (N')^{1-2\beta}}{\varepsilon^2} + \frac{1}{\varepsilon^2 (N')^{2\beta+1}} \right),$$

and the choice of $\delta \asymp \varepsilon^{2/(2\beta+1)}$ yields the boundedness of the last expression and hence the desired result. $\square$

The proofs of Theorems 3 and 4 follow the same steps as the proofs of Theorems 1 and 2 with slight modifications and are omitted.

A. GOLDENSHLUGER
DEPARTMENT OF STATISTICS
UNIVERSITY OF HAIFA
HAIFA 31905
ISRAEL
E-MAIL: goldensh@stat.haifa.ac.il

A. TSYBAKOV
LABORATOIRE DE PROBABILITÉS
ET MODÈLES ALÉATOIRES
UNIVERSITÉ PARIS VI
4 PLACE JUSSIEU
PARIS 75252
FRANCE
E-MAIL: tsybakov@ccr.jussieu.fr

A. ZEEVI
GRADUATE SCHOOL OF BUSINESS
COLUMBIA UNIVERSITY
3022 BROADWAY
NEW YORK 10027
USA
E-MAIL: assaf@gsb.columbia.edu